\documentclass[a4paper,12pt]{article}
\usepackage{amsfonts,amssymb, pdfpages,xcolor}
\usepackage[latin1]{inputenc}
\usepackage[normalem]{ulem}

\title{Approximate  solutions to the initial value problem for some  compressible flows in presence of shocks and void regions. }
\author{M. Colombeau*,\\ \textit{mcolombeau@ime.usp.br}\\
 Instituto de Matematica e Estatistica,\\Universidade  de S\~ao Paulo, Brazil.}

\date{}
\begin{document}
\maketitle

\begin {abstract} For the natural initial conditions $L^1$ in the density field (more generally a positive bounded Radon measure) and  $L^\infty$ in the  velocity field we obtain global approximate solutions to the Cauchy problem  for the 3-D systems of isothermal and isentropic gases,  the 2-D shallow water equations and the 3-D system of collisionnal self-gravitating gases. We obtain a sequence of  functions which are differentiable in time and continuous in space and tend to satisfy the equations in the sense of distributions in the space variables and in the strong sense in the time variable.  The method of construction  relies on the study of a specific family of two ODEs in a classical Banach space (one for the continuity equation and one for the Euler equation). Standard convergent numerical methods for the solution of these ODEs can be used to  provide concrete approximate solutions. It has been checked in numerous cases in which the solutions of systems of fluid dynamics are known that our constuction   always gives back the known solutions. It is also proved it gives the classical analytic solutions in the domain of application of the Cauchy-Kovalevska theorem.

\end{abstract}

\textit{* This research has been  done thanks to  financial support of FAPESP, processo 2012/15780-9.}\\

\textbf{1.  Introduction}.\\

 Continuying a previous study of self gravitating pressureless fluids \cite{ColombeauODE} we construct sequences of  approximate solutions   for 3-D compressible  isothermal gases and isentropic gases in presence of shocks and void regions; this applies also to  the shallow water equations and to the system of collisionnal selfgravitating isothermal gases. The need of the search of mathematical solutions to the initial value problem for compressible fluids in several dimensions when shocks show up has been recently pointed out by P. D. Lax  and D. Serre \cite{Lax1, Serre}. \\

The starting point of our attempt  \cite{ColombeauSiam, ColombeauNMPDE, Colombeauideal} was the search of a numerical scheme permitting mathematical proofs  of partial results aiming to connect the functions produced by  the scheme and the system one tries to solve: indeed, in absence of well defined solutions, it is natural to try to prove that the  "approximate solutions" from the scheme tend to satisfy the equations in the sense of distributions, which makes sense in absence of properly defined mathematical solution, see the appendix.  Sequences of functions that tend to satisfy the equations could be a provisional substitute of solutions, since one has to cope with an  absence of known usual weak solutions  in 2-D and 3-D in presence of shocks \cite{Lax1} p. 143, \cite{Serre} p. 143 and \cite{Novotny} p. 150. These sequences  permit in particular to explain results observed from computing \cite{Lax1, Lax2}. In order to obtain full proofs, we replaced the original scheme  \cite{ColombeauSiam, ColombeauNMPDE, Colombeauideal} by ODEs in Banach spaces whose solutions provide the sequence of "approximate solutions", at the same time as one retains the possibility of computing the solutions of these ODEs to check that the method gives  the correct results.\\

From a physical viewpoint the equations of fluid dynamics are mared with some imprecision since they do not take into account some minor effects and the molecular structure of matter. It is natural to expect these equations and their imprecision should be stated in the sense of distributions in the space variables. Sequences $[(x,y,z,t)\longmapsto U(x,y,z,t,\epsilon)]_\epsilon$ of approximate solutions in the sense of distributions   enter into this imprecision for $\epsilon>0$ small enough. Therefore  these sequences appear to provide a natural representation of physical solutions.
Of course one has checked that the sequences of approximate solutions we construct always give the correct solution in all the numerous tests that were performed. Uniqueness of the limit when $\epsilon\rightarrow 0$ for a general class of sequences of approximate solutions containing those we construct has been proved  in the  case of linear systems with analytic coefficients \cite{ColombeauCK}. For the equations of fluid dynamics considered here, in absence of a known similar uniqueness result in a suitable family of approximations, one has to content with numerical confirmations.\\

These sequences of approximate solutions have been introduced under the name of  weak asymptotic solutions   by Danilov, Omel'yanov and Shelkovich \cite{Danilov1}, as an extension of Maslov asymptotic analysis, and they have proved to be an efficient mathematical tool to study creation and superposition of singular solutions  to various nonlinear PDEs \cite {Danilov1,Mitrovic, Shelkovich2, Shelkovich3, Panov, ShelkovichRMS,Shelkovichmat,Shelkovich1}. A weak asymptotic solution for the system
 $$\frac{\partial}{\partial t}U_i+\frac{\partial}{\partial x}f_i(U_1,\dots, U_p)+\frac{\partial}{\partial y}g_i(U_1,\dots, U_p)=0, 1\leq i\leq p,\ U_i:\mathbb{R}^2\times\mathbb{R}^+\mapsto\mathbb{R},$$ is a sequence $(U_{1,\epsilon},\dots, U_{p,\epsilon})_\epsilon$ of functions such that $\forall i\in\{1,\dots,p\},\  \forall \psi\in \mathcal{C}_c^\infty(\mathbb{R}^2)  \  \forall t $
 \begin{equation} \ \ \int_{\mathbb{R}^2}[\frac{\partial}{\partial t}U_{i,\epsilon}.\psi-f_i(U_{1,\epsilon},\dots, U_{p,\epsilon})\frac{\partial}{\partial x}\psi-g_i(U_{1,\epsilon},\dots, U_{p,\epsilon})\frac{\partial}{\partial y}\psi]dxdy\rightarrow 0\end{equation} 
when $\epsilon\rightarrow 0$ 
 if we consider the strong derivative in $t$ and weak derivatives in $x,y$ as this will be the case in this paper. In short the sequence $(U_{1,\epsilon},\dots, U_{p,\epsilon})_\epsilon$ tends to satisfy the equations. Of course the $U_{i,\epsilon}'s$ are also chosen so as to satisfy the initial condition stated in a natural sense.\\

 The approximate solutions that we construct satisfy uniform $L^1$ bounds in density (consequence of mass conservation) but, in presence of pressure, one does not obtain a bounded  velocity; this should not be too much unexpected: indeed  various instances are known in which the ideal inviscid equations give infinite velocity  \cite{Guyon} sections 6.6.2 , 7.4.3, \cite{Campos} section 14.9.1. The concept of sequences of approximate solutions in this paper permits infinite velocity at the limit $\epsilon\rightarrow 0$. \\

 We construct a weak asymptotic method for the 3-D systems of  isothermal and isentropic gases that we state in 2-D for convenience (immediate extension to 3-D) in the form 
 \begin{equation}\frac{\partial}{\partial t}\rho+\frac{\partial}{\partial x}(\rho u)+\frac{\partial}{\partial y}(\rho v)=0,   \end{equation}
  \begin{equation}  \frac{\partial}{\partial t}(\rho u)+\frac{\partial}{\partial x}(\rho u^2)+\frac{\partial}{\partial y}(\rho uv)+\frac{\partial}{\partial x}p=0,  \end{equation}
  \begin{equation}\frac{\partial}{\partial t} (\rho v)+\frac{\partial}{\partial x}(\rho uv)+\frac{\partial}{\partial y}(\rho v^2)+\frac{\partial}{\partial y}p=0,  \end{equation}
  \begin{equation} p=K \rho, \ \ (isothermal \  gases),  \end{equation}
 \begin{equation} p=K \rho^\gamma, \ \ 1<\gamma\leq 2, \ \ (isentropic \ gases).  \end{equation}
   
   The notation is: $\rho=\rho(x,y,t)$ is the gas density, $(u,v)=(u(x,y,t),v(x,y,t))$ is the velocity vector in the $x,y$ directions respectively, $p=p(x,y,t)$ is the pressure and $K$ is a constant.  The problem is studied on the $n$-dimensional torus $\mathbb{T}^n=\mathbb{R}^n/(2\pi\mathbb{Z})^n, \ n=1,2,3.$\\

We combine the system of isothermal  gases and self-gravitation to obtain the classical system of collisionnal selfgravitating gases which models Jeans gravitational instability \cite{Charru, Coles}: a large enough cloud of gas possibly at rest collapses gravitationally besides  pressure untill an equilibrium is reached with pressure, forming a star or a planet. A mathematical proof in this paper asserts that the numerical result observed is a depiction of an approximate solution of the equations with arbitrary precision, only limited by the cost of calculations.\\

   Then, since the equations are similar, we construct  weak asymptotic solutions for the 2-D shallow water equations stated in the form
 \begin{equation}\frac{\partial}{\partial t}h+\frac{\partial}{\partial x}(hu)+\frac{\partial}{\partial y}(hv)=0,   \end{equation}
  \begin{equation}\frac{\partial}{\partial t} (hu)+\frac{\partial}{\partial x}(hu^2)+\frac{\partial}{\partial y}(huv)+gh\frac{\partial}{\partial x}(h+a)=0,  \end{equation}
  \begin{equation} \frac{\partial}{\partial t}(hv)+\frac{\partial}{\partial x}(huv)+\frac{\partial}{\partial y}(hv^2)+gh\frac{\partial}{\partial y}(h+a)=0,  \end{equation}
 \\
 where $h=h(x,y,t)$ is the water elevation, $(u,v)=(u(x,y,t),v(x,y,t))$ is the velocity vector in the $x,y$ directions respectively, $a=a(x,y)$ is the bottom elevation assumed to be of class $\mathcal{C}^2$ and $g=9.8$. The problem is studied on the $n$-dimensional torus $\mathbb{T}^n, \ n=1,2.$\\

  The method of proof is the statement and study of a system of ODEs (one for the continuity equation, one for the Euler equation) in the Banach space $\mathcal{C}(\mathbb{T}^n), \ n=1,2,3$, of  continuous functions. We expose the proof in 1-D since the extension to several dimensions is straightforward. One uses the theory of ODEs in the Lipschitz case and a priori estimates needed to prove existence of a global in positive time solution. It has been inspired by the method introduced in \cite{ColombeauODE} in the pressureless case and in \cite{ColombeauJDE} to put in evidence continuations of the analytic solutions after the time of the analytic blow up.\\ 
   
  Numerical calculations  of these approximate solutions can be done easily by solving the above ODEs in Banach space by classical convergent numerical methods for ODEs such as the explicit Euler order one method or the  Runge Kutta RK4 method. This has  given the known solutions in all tests, which could have been expected since the method in this paper is issued from the numerical scheme in \cite{ColombeauSiam,ColombeauNMPDE,Colombeauideal}, for which a large amount of numerical verifications for initial conditions in which the solution is known  has been done on classical and demanding tests of different nature (Sod\cite{Sod}, Woodward-Colella\cite {WoodwardColella}, Toro\cite{Toro, Toro2}, Lax \cite{Lax1,Lax2}, Bouchut-Jin-Li\cite{Bouchut}, LeVeque\cite{LeVeque}, Cherkov-Kurganov-Rykov\cite{Chertock}, Evje-Flatten\cite{EvjeFlatten}).\\

  
 \textbf{2. Approximate solutions of  the system of isothermal gases and numerical confirmations.} \\
 
  We state the classical system of isothermal gases  in 1-D in the form:

\begin{equation}\frac{\partial}{\partial t}\rho+\frac{\partial}{\partial x}(\rho u)=0,   \end{equation}
  \begin{equation}\frac{\partial}{\partial t} (\rho u)+\frac{\partial}{\partial x}(\rho u^2)+\rho\frac{\partial}{\partial x}\Phi=0,  \end{equation}
   \begin{equation} \Phi=K log(\rho), \ K\geq 0  \ a \ given \  constant.  \end{equation}
   Therefore the last term in the Euler equation (11) is in the familiar form $\rho \frac{\partial}{\partial x}\Phi=\frac{\partial}{\partial x}p,$ with   $ p=K\rho$.  Here $\Phi$ is the density of the body force caused by the pressure: $\rho \vec{\nabla}\Phi=\vec{\nabla}p$.  Formulation (12) excludes void regions. This apparent defect will be repaired by the term $\rho\vec{\nabla}\Phi$  in the Euler equation: it will be proved theoretically and checked numerically that the construction works even in presence of void regions. We use the notation
  
   
   \begin{equation} u^+=max(u,0), \ u^-=max(-u,0),\end{equation}
 so that
 \begin{equation} u=u^+-u^-, \  |u|=u^++u^-. \end{equation}
  We study the system in 1-D since the extensions to 2-D and 3-D are straightforward, following these extensions in section 6 of \cite{ColombeauODE}. The 2-D extensions of formulas (15, 16) below are given as (44, 45) below.
 In 1-D we state the  method as the system of ODEs (15, 16) complemented by the formulas (17, 18):
  \begin{equation}\frac{d}{dt}\rho(x,t,\epsilon)=\frac{1}{\epsilon}[(\rho u^+)(x-\epsilon,t,\epsilon)-(\rho |u|)(x,t,\epsilon)+(\rho u^-)(x+\epsilon,t,\epsilon)]+ \epsilon^\beta,\end{equation}
for some $\beta>0$ to be made precise later,
\begin{equation}\frac{d}{dt}(\rho u)(x,t,\epsilon)=\frac{1}{\epsilon}[(\rho u u^+)(x-\epsilon,t,\epsilon)-(\rho u |u|)(x,t,\epsilon)+(\rho u  u^-)(x+\epsilon,t,\epsilon)] -\rho(x,t,\epsilon)\frac{\partial}{\partial x}\Phi(x,t,\epsilon),\end{equation}
\begin{equation} u(x,t,\epsilon)=\frac{(\rho u)(x,t,\epsilon)}{\rho(x,t,\epsilon)},\end{equation}
for which we will prove that $\rho(x,t,\epsilon)>0,$ thus permitting division. The term $\epsilon^\beta$ in (15) is needed because of  the specific form of the state law of isothermal gases and is not needed for the  isentropic gases and shallow water equations.\\

 We approximate the state law (12) in the form
  \begin{equation}\Phi(x,t,\epsilon)=K[log(\rho(.,t,\epsilon)+\epsilon^N)*\phi_{\epsilon^\alpha}](x), \ N>0, \end{equation} 
  in which the term $\epsilon^N$  has been introduced to permit  void regions. 
  We introduce an auxiliary function  $\phi\in \mathcal{C}_c^\infty(\mathbb{R}), \phi\geq 0, \phi(-x)=\phi(x) \ \forall x$ and $ \int\phi(x)dx=1.$ If $\mu>0$ we set $\phi_{\mu}(x)=\frac{1}{\mu}\phi(\frac{x}{\mu}),$ so that the family $\{\phi_{\mu}\}_{\mu}$ tends in the sense of distributions to the Dirac measure when $\mu\rightarrow 0^+$. We use a real number $\alpha, 0<\alpha<1$, to be made more precise later. The convolution in (18) is justified by the fact (12) is a state law and, as such,  is physically valid only in space regions  larger than those in which the basic conservation laws (10, 11) are valid.\\
 

 We first  establish a priori inequalities to prove existence of a global solution to (15-18). For fixed $\epsilon>0$ we assume existence of a solution
$$[0,\delta(\epsilon)[ \longmapsto  (\mathcal{C}(\mathbb{T}))^2$$
\begin{equation} \ \  \ \ \  \ \ \ \ \ \ \ \ \ \ \   t  \longmapsto[x\mapsto(\rho(x,t,\epsilon),(\rho u)(x,t,\epsilon))]\end{equation}
continuously differentiable  on $[0,\delta(\epsilon)[$ (with a right hand-side derivative at $t=0$) having the following properties for each fixed $\epsilon$: 
  \begin{equation} \exists m>0 \ / \ \rho(x,t,\epsilon)\geq m \ \forall x\in \mathbb{T} \ \forall t\in [0,\delta(\epsilon)[,  \end{equation}
  \begin{equation}  \exists M>0 \ / \ \|u(.,t,\epsilon)\|_\infty \leq M, \ \ \|\rho(.,t,\epsilon)\|_\infty \leq M \ \forall t\in [0,\delta(\epsilon)[.
  \end{equation}
  Note that $m$ and $M$  depend on $\epsilon$.\\

    The aim of the following a priori inequalities is to obtain bounds independent on $m$ and $M$ to replace (20, 21) for fixed $\epsilon>0$, depending   on the initial condition and on $\epsilon$,  in order to prove that the solution can be extended for $t>\delta(\epsilon)$. All constants, denoted $const$ in the proposition below, depend on the initial condition and possibly on a bounded time interval,  but not on the values $m$ and $M$ in the a priori assumption (20, 21) and not on $\epsilon$. The independence on $\epsilon$  will be basic when they will  be used later at the limit $\epsilon\rightarrow 0$ to prove that the solutions of the ODEs tend to satisfy the equations. The dependence of $const$ in $t$ is explicitely stated when this appears clearer.\\ 
  
  We assume $\rho_0$ and $u_0$ are given initial conditions with the properties $\rho_0 \in L^1(\mathbb{T})$ and $ u_0\in   L^\infty(\mathbb{T})$ and that $\rho_0^\epsilon$ and $ u_0^\epsilon$ are regularizations of $\rho_0$ and  $u_0 \in \mathcal{C}(\mathbb{T})$, with  uniform $L^1$ and $L^\infty$ bounds respectively (independent on $\epsilon$), and $\rho_0^\epsilon (x)>0 \  \  \forall x\in \mathbb{T}$.\\
  
  \textbf {Proposition 1 (a priori inequalities).}\\
    \begin{equation} \bullet \forall       t\in [0,\delta(\epsilon)[ \  \ \int_{-\pi}^\pi\rho(x,t,\epsilon)dx = \int_{-\pi}^\pi\rho_0^\epsilon(x)dx+2\pi\epsilon^\beta t, \ \ \ \ \ \ \ \ \ \ \ \ \ \ \ \ \ \ \ \ \ \ \ \end{equation}  
     \begin{equation}  \bullet\forall       t\in [0,\delta(\epsilon)[ \ \ \|\frac{\partial}{\partial x}\Phi(.,t,\epsilon)\|_\infty \leq \frac{const}{\epsilon^{3\alpha}},  \ \ \ \ \ \ \ \ \ \ \ \ \ \ \ \ \ \ \ \ \ \ \  \ \ \ \ \ \ \ \ \ \ \ \ \ \ \ \ \ \ \ \ \ \ \ \end{equation}
      \begin{equation} \bullet\forall       t\in [0,\delta(\epsilon)[ \ \ \|u(.,t,\epsilon)\|_\infty \leq  \|u_0^\epsilon\|_\infty +\frac{const} {\epsilon^{3\alpha}}t. \ \ \ \ \ \ \ \ \ \ \ \ \ \ \ \ \ \ \ \ \ \ \  \ \ \ \ \ \ \ \ \ \ \  \end{equation}
       Set 
  \begin{equation} k(\epsilon)=\|u_0^\epsilon\|_\infty+  \frac{const. \delta(\epsilon)}{\epsilon^{3\alpha}}.\end{equation}
   Then   
      \begin{equation} \bullet  \forall       t\in [0,\delta(\epsilon)[, \forall x\in \mathbb{R}, \ \rho_0^\epsilon(x)
exp(\frac{-k(\epsilon)t}{\epsilon}) \leq \rho(x,t,\epsilon) \leq 2\|\rho_0^\epsilon\|_\infty exp(\frac{2k(\epsilon)t}{\epsilon}).\end{equation} 
  \\  
 \textit{proof}. From (14, 15)\\

      $\frac{d}{dt}\int_{-\pi}^{+\pi}\rho(x,t,\epsilon)dx=
      \frac{1}{\epsilon}[\int_{-\pi}^{+\pi}(\rho u^+)(x-\epsilon,t,\epsilon)dx-\int_{-\pi}^{+\pi}(\rho u^+)(x,t,\epsilon)dx-\int_{-\pi}^{+\pi}(\rho u^-)(x,t,\epsilon)dx\\
\\
+\int_{-\pi}^{+\pi}(\rho u^-)(x+\epsilon,t,\epsilon)dx]+2\pi\epsilon^\beta 
   =0+2\pi\epsilon^\beta$\\ 
 \\    
     by periodicity of $\rho$ and $u$. \\ 
     
    Inequality  (23) is proved as follows:  (18) implies
  $$ \frac{\partial}{\partial x}\Phi(x,t,\epsilon)=K\int log[\rho(x-y,t,\epsilon)+\epsilon^N]\frac{1}{\epsilon^{2\alpha}}  \phi'(\frac{y}{\epsilon^\alpha})dy.$$
  Using  (22) when $\rho(x-y,t,\epsilon)\geq 1$ (then $log(\rho+\epsilon^N)<\rho$), and, using that  the above $log$ is bounded in absolute value by $const.log(\frac{1}{\epsilon^N})$  when  
 $\rho(x-y,t,\epsilon)< 1$, one obtains
 $$|\frac{\partial}{\partial x}\Phi(x,t,\epsilon)|\leq const.  log(\frac{1}{\epsilon}) \frac{1}{\epsilon^{2\alpha}}  \|\phi'\|_\infty \leq const\frac{1}{\epsilon^{3\alpha}},$$
  which gives the desired result ($const$ does not depend on $t$ from (22) when $t$ ranges in a bounded interval). \\
 
  Now let us prove inequality (24). From (15) and the assumption that the solution of the ODE  is of class $\mathcal{C}^1$ on $[0,\delta(\epsilon)[$, valued in the Banach space $\mathcal{C}(\mathbb{T})$, one obtains for fixed $\epsilon>0$ and for $dt>0$ small enough ($t+dt<\delta(\epsilon))$ that\\

     $\rho(x,t+dt,\epsilon)=\rho(x,t,\epsilon)+$\\
     $$\frac{dt}{\epsilon}[(\rho u^+)(x-\epsilon,t,\epsilon)-(\rho |u|)(x,t,\epsilon)+(\rho u^-)(x+\epsilon,t,\epsilon)]+dt.o(x,t,\epsilon)(dt)+\epsilon^\beta dt=$$ \begin{equation}\frac{dt}{\epsilon}(\rho u^+)(x-\epsilon,t,\epsilon)+(1-\frac{dt}{\epsilon}|u|(x,t,\epsilon))\rho(x,t,\epsilon)+\frac{dt}{\epsilon}(\rho u^-)(x+\epsilon,t,\epsilon)+dt.o(x,t,\epsilon)(dt)+\epsilon^\beta dt\end{equation}
     \\ 
      where $\|o(.,t,\epsilon)(dt)\|_\infty \rightarrow 0$ when $dt\rightarrow 0$ uniformly for $t$ in a compact set of $[0,\delta(\epsilon)[$, from the mean value theorem under the form $\|f(t+dt) -f(t)-f'(t)dt\|\leq sup_{0<\theta<1}\|f'(t+\theta dt)-f'(t)\||dt|$. Notice that there is no uniformness in $\epsilon$. For $dt>0$ small enough (depending on $\epsilon$)   the single term $  (1-\frac{dt}{\epsilon}|u(x,t,\epsilon)|)\rho(x,t,\epsilon)$ dominates the term $dt.o(x,t,\epsilon)(dt)$ from (20, 21). Therefore since $\rho u^\pm\geq0$ one can invert (27):\\
 
$\frac{1}{\rho(x,t+dt,\epsilon)}=$
$[\frac{dt}{\epsilon}(\rho u^+)(x-\epsilon,t,\epsilon)+[1-\frac{dt}{\epsilon}|u|(x,t,\epsilon)]\rho(x,t,\epsilon)+\\
\\
\frac{dt}{\epsilon}(\rho u^-)(x+\epsilon,t,\epsilon)+\epsilon^\beta dt]^{-1}+dt.o(x,t,\epsilon)(dt)$\\
\\
where the new $o$ has still the property that $\|o(.,t,\epsilon)(dt)\|_\infty \rightarrow 0$ when $dt\rightarrow 0$ uniformly for $t\in[0,\delta'] $  if $\delta'<\delta(\epsilon)$.\\
      
      
      Applying the analog of (27) with $\rho u$ in place of $\rho$, with the supplementary term $\rho \frac{\partial}{\partial x}\Phi$ from (16)  and absence of the term $\epsilon^\beta$,  one obtains (from 20, 21):
 $$u(x,t+dt,\epsilon)=\frac{(\rho u)(x,t+dt,\epsilon)}{\rho(x,t+dt,\epsilon)}=$$
$$ \frac
 {\frac{dt}{\epsilon}(\rho u u^+)(x-\epsilon,t,\epsilon)+[1-\frac{dt}{\epsilon}|u|(x,t,\epsilon)](\rho u)(x,t,\epsilon)+\frac{dt}{\epsilon}(\rho u u^-)(x+\epsilon,t,\epsilon)}
 {\frac{dt}{\epsilon}(\rho u^+)(x-\epsilon,t,\epsilon)+[1-\frac{dt}{\epsilon}|u|(x,t,\epsilon)]\rho(x,t,\epsilon)+\frac{dt}{\epsilon}(\rho u^-)(x+\epsilon,t,\epsilon)+\epsilon^\beta dt}$$
 \begin{equation} -dt\frac{\rho(x,t,\epsilon)}{\rho(x,t+dt,\epsilon)}\frac{\partial}{\partial x}\Phi(x,t,\epsilon)+ dt.o(x,t,\epsilon)(dt)\end{equation}
where the new $o$ has the same property as in (27)  for fixed $\epsilon$.
 For $dt>0$ small enough the first term  in the second member is  a  barycentric combination of $u(x-\epsilon,t,\epsilon),  u(x,t,\epsilon) $ and $ u(x+\epsilon,t,\epsilon)$ (which are in numerator in factor of $\rho$ inside $\rho u$), dropping the term $\epsilon^\beta dt$. For fixed $\epsilon$ the quotient  $\frac{\rho(x,t+dt,\epsilon)}{\rho(x,t,\epsilon)} $ tends to $1$ when $dt\rightarrow 0$ (use (20, 21 and 27). Finally it follows from (23) and (28) that 
 \begin{equation} \|u(.,t+dt,\epsilon)\|_\infty\leq  \|u(.,t,\epsilon)\|_\infty+dt\frac{const}{\epsilon^{3\alpha}}+dt.\|o(.,t,\epsilon)(dt)\|_\infty\end{equation}
 with uniform bound of $o$ when $t$ ranges in a compact set in $[0,\delta(\epsilon)[$. One obtains the bound (24)  by dividing the interval $[0,t]$ into $n$ intervals $[\frac{it}{n},\frac{(i+1)t}{n}], 0\leq i \leq n-1$, and applying (29) in each subinterval:
 $$\|u(.,(i+1)\frac{t}{n},\epsilon)\|_\infty \leq \|u(.,i\frac{t}{n},\epsilon)\|_\infty+\frac{t}{n}\frac{const}{\epsilon^{3\alpha}}+\frac{t}{n}o(\frac{t}{n}),$$
 summing on $i$ and using that $o(\frac{t}{n})\rightarrow 0$ when $n\rightarrow \infty$, see more details in \cite{ColombeauODE}.
  \\


 Now let us prove inequalities (26). From (15), $$\frac{d}{dt}\rho(x,t,\epsilon)\geq -\frac{1}{\epsilon}(\rho |u|)(x,t,\epsilon),$$ since $\rho,u^+$ and $u^-$ are positive (13, 20). Therefore, from (24, 25),
\begin{equation}\frac{d}{dt}\rho(x,t,\epsilon)\geq -\frac{k(\epsilon)}{\epsilon}\rho(x,t,\epsilon).\end{equation}
Let $v(x,t,\epsilon)=\rho(x,0,\epsilon)exp(-\frac{k(\epsilon)}{\epsilon}t)$. Then, using assumption (20) to divide by $\rho$,\\
\begin{equation}\frac{\frac{d}{d t}\rho}{\rho}(x,t,\epsilon)\geq\frac{\frac{d}{d t}v}{v}(x,t,\epsilon)=-\frac{k(\epsilon)}{\epsilon}.\end{equation}
 By integration, since $\rho$ and $v$ have same initial condition and are positive, $log(\rho)\geq log (v)$, i.e. $\rho(x,t,\epsilon)\geq v(x,t,\epsilon)$, i.e.
\begin{equation}\rho(x,t,\epsilon)\geq \rho(x,0,\epsilon)exp(-\frac{k(\epsilon)}{\epsilon} t),\end{equation}
which is the left hand-side  inequality  (26). Now let us prove the right hand-side  inequality  (26).\\

From the positiveness of the two terms $\rho u^{\pm}$ in (15) 

\begin{equation} \rho(x,t,\epsilon)\leq \rho_0(x,\epsilon)+\frac{2}{\epsilon}\int_0^t\|\rho(.,s,\epsilon)\|_\infty \|u(.,s,\epsilon)\|_\infty ds.   \end{equation}
From (24, 25), if $t\in [0,\delta(\epsilon)[$
$$\rho(x,t,\epsilon)\leq \| \rho_0(.,\epsilon)\|_\infty+\frac{2}{\epsilon}\int_0^t\|\rho(.,s,\epsilon)\|_\infty k(\epsilon) ds.$$

 Since this holds for all $x$
\begin{equation} \|\rho(.,t,\epsilon)\|_\infty\leq \| \rho_0(.,\epsilon)\|_\infty+\frac{2}{\epsilon}k(\epsilon)\int_0^t\|\rho(.,s,\epsilon)\|_\infty ds.   \end{equation}
Gronwall's inequality implies
\begin{equation} \|\rho(.,t,\epsilon)\|_\infty\leq \| \rho_0(.,\epsilon)\|_\infty exp(\frac{2}{\epsilon}k(\epsilon) t).\end{equation}
$\Box$ \\   
 
 For fixed $\epsilon>0$, if $0<\lambda<1$ and $\Omega_\lambda:=\{(X,Y)\in \mathcal{C}(\mathbb{T})^2 / \forall x\in \mathbb{R} \ \lambda<X(x)<\frac{1}{\lambda}, |Y(x)|<\frac{1}{\lambda}\}$ the equations (15-18) with variables $X=\rho, Y=\rho u$ have the Lipschitz property on $\Omega_\lambda$ with values in $\mathcal{C}(\mathbb{T})$, with Lipschitz constants uniform in $\Omega_\lambda$, see \cite {ColombeauODE} section 4 for details. The existence of a unique global solution to (15-18) for fixed $\epsilon$ is obtained from the a priori estimates in proposition 1 from classical  arguments of the theory of ODEs in Banach spaces in the Lipschitz case as exposed in section 4 of \cite{ColombeauODE}.\\
 
  It remains to prove that the solution of the system of ODEs (15, 16) complemented by (17, 18) provides a weak asymptotic method for system (10-12) when $\epsilon\rightarrow 0$. To this end one has to prove that $\forall t>0, \forall \psi\in \mathcal{C}_c^\infty(\mathbb{R})$ (36-38) below hold when $\epsilon\rightarrow 0$:
 \begin{equation} \int \frac{d}{dt}\rho(x,t,\epsilon)\psi(x)dx=\int (\rho u)(x,t,\epsilon)\psi'(x)dx +f(\epsilon),\end{equation}
  \begin{equation} \int \frac{d}{dt}(\rho u)(x,t,\epsilon)\psi(x)dx=\int (\rho u^2)(x,t,\epsilon)\psi'(x)dx -\int \rho(x,t,\epsilon)\frac{\partial}{\partial x}\Phi(x,t,\epsilon)\psi(x)dx+f(\epsilon),\end{equation}
  \begin{equation}  \int\Phi(x,t,\epsilon)\psi(x)dx=\int K \log[\rho(x,t,\epsilon)]\psi(x)dx+f(\epsilon), \end{equation} 
 where the three different $f(\epsilon)$ tend to 0 when $\epsilon\rightarrow 0$.\\

    The proof of (36) is as follows: from (14, 15, 22, 24), a change of variable and $\frac{\psi(x+\epsilon)-\psi(x)}{\epsilon}=\psi'(x)+O_x(\epsilon)$,\\ 
     
  $ \int \frac{d}{dt}\rho(x,t,\epsilon)\psi(x)dx=
  \frac{1}{\epsilon}\int (\rho u^+)(x,t,\epsilon)[\psi(x+\epsilon)-\psi(x)]dx -\frac{1}{\epsilon}\int (\rho u^-)(x,t,\epsilon)$\\
\\
$[\psi(x)-\psi(x-\epsilon)]dx+\int\epsilon^\beta \psi(x)dx =
  \int (\rho u)(x,t,\epsilon)\psi'(x)dx+\int_{compact} (\rho u^+)(x,t,\epsilon) O_x(\epsilon)dx+$\\
\\
$ \int_{compact} (\rho u^-)(x,t,\epsilon) O_x(\epsilon)dx+O(\epsilon^\beta)= \int(\rho u)(x,t,\epsilon)\psi'(x) dx+(const+\frac{const}{\epsilon^{3\alpha}}t) O(\epsilon)$\\
\\
$+O(\epsilon^\beta)= \int(\rho u)(x,t,\epsilon)\psi'(x) dx+O(\epsilon^{1-3\alpha})+O(\epsilon^\beta)$.\\
  \\
   This gives (36) if $0<\alpha<\frac{1}{3}$. The proof of (37) is similar since the additional terms $\rho(x,t,\epsilon)\frac{\partial}{\partial x}\Phi(x,t,\epsilon)$   in (16) and (37) simplify. One obtains a remainder $\frac{const}{\epsilon^{6\alpha}} t O(\epsilon)$ because of one more factor $u$ and the bound (24). Finally one chooses $0<\alpha<\frac{1}{6}$.\\

 To check (38) one has to prove from (18) that 
  \begin{equation}\forall \psi \in \mathcal{C}_c^\infty(\mathbb{R}) \ \int\{[log(\rho(.,t,\epsilon)+\epsilon^N)*\phi_{\epsilon^\alpha}](x)-log[\rho(x,t,\epsilon)]\}\psi(x)dx \rightarrow 0   \end{equation}
  when $\epsilon\rightarrow 0$. To this end we share the integral (39) into the two parts (40, 41) below and we prove that each of them tends to 0 when $\epsilon\rightarrow 0$. Let 
   \begin{equation}I=  \int\{[log(\rho(.,t,\epsilon)+\epsilon^N)*\phi_{\epsilon^\alpha}](x)-log[\rho(x,t,\epsilon)+\epsilon^N]\}\psi(x)dx \end{equation} 
  and 
  \begin{equation} J=\int\{log[\rho(x,t,\epsilon)+\epsilon^N]-log[\rho(x,t,\epsilon)]\}\psi(x)dx.\end{equation}
  Now\\

  $I=\int\{log[\rho(x-\epsilon^\alpha\mu,t,\epsilon)+\epsilon^N]-log[\rho(x,t,\epsilon)+\epsilon^N]\}\phi(\mu)\psi(x)d\mu dx=$\\
\\
  $\int log[\rho(x,t,\epsilon)+\epsilon^N]\phi(\mu)[\psi(x+\epsilon^\alpha\mu)-\psi(x)]d\mu dx.$\\
  \\
  Since $\rho(x,t,\epsilon)\geq 0$ from (26), using (22) in the case $\rho(x,t,\epsilon)>1$ and using $\epsilon^N$ in the case  $\rho(x,t,\epsilon)\leq1$, as in the proof of (23),  one has  $|I|\leq const.log(\frac{1}{\epsilon})\epsilon^\alpha$; therefore $I\rightarrow 0$ when $\epsilon\rightarrow 0$.\\
 
  Now (41) and the mean value theorem give 
  \begin{equation}|J|\leq \epsilon^N\frac{1}{min (\rho)} const\end{equation}
  if $min(\rho)$ denotes the inf of $\rho(x,t,\epsilon)$ for fixed  $t$ and $\epsilon$ when $x$ ranges in $\mathbb{T}$.
  The problem is to obtain a suitable inf. bound of $min(\rho)$. The term $\epsilon^\beta$ in (15) has been introduced for this purpose. Indeed, from (15),
  $$\frac{d\rho}{dt}(x,t,\epsilon)\geq-\frac{1}{\epsilon}\rho(x,t,\epsilon)\|u(.,t,\epsilon)\|_\infty +\epsilon^\beta,$$
  i.e. setting $A:=const \frac{1}{\epsilon^{1+3\alpha}}$ and  $\ B:=\epsilon^\beta$, the fact from (24) that  $\|u\|_\infty \leq \frac{const}{\epsilon^{3\alpha}}$ on any bounded time interval implies that on such interval 
  $$\frac{d\rho}{dt}\geq -A\rho+B.$$
   If one considers the ODE $\frac{dX}{dt}(x,t)=-AX(x,t)+B$ with initial condition $X(x,0)=\rho_0^\epsilon(x)$, its solution is
   \begin{equation}X(x,t)=\rho_0^\epsilon(x) e^{-At}+\frac{B}{A}(1-e^{-At}).\end {equation}
   For fixed $ t>0$ and for $\epsilon>0$ small enough
   $$X(x,t)\geq \frac{B}{2A}=const(t).\epsilon^{1+\beta+3\alpha},$$
   in which we state the constant as $const(t)$ since this constant depends on $t$ through (24). We have
  $$\rho(x,t,\epsilon)\geq X(x,t),$$
  therefore
  $$\rho(x,t,\epsilon)\geq const(t).\epsilon^{1+\beta+3\alpha}.$$
\\
  From (42)
  $|J|\leq const(t).\epsilon^{N-1-\beta-3\alpha}$ and  it suffices to choose $\beta+3\alpha<N-1$ 
  to have that $J\rightarrow 0$ when $\epsilon \rightarrow 0.$\\
  
   It follows from the proof that the  possible presence of void regions, approximated in the initial conditions by  $\rho_0^\epsilon(x)\geq \epsilon  \ \forall x$, does not cause a problem although they appear a priori excluded by the formulation (12). Note also that the presence of $\epsilon$ in a denominator in (26) allows concentrations of matter, that can also be accepted in initial conditions by choosing $\rho_0^\epsilon(x)\leq \frac{1}{\epsilon}  \ \forall x$:  the initial condition $\rho_0$ can even be a bounded Radon measure. One could notice that  formula (24) allows the possibility of infinite velocity at the limit $\epsilon=0$. This should not be troublesome: the remark that the ideal equations (2-6) could lead to solutions with infinite velocity in certain circumstances has been known since long time \cite{Campos} section 14.9.1, \cite{Guyon} section 6.6.2, 7.4.3.\\

   Finally we have proved that under the initial conditions $\rho_0\in L^1 (\mathbb{T}),  u^0\in L^\infty (\mathbb{T}), \  \rho_0\geq 0$,  and approximating the initial conditions by a family $(\rho_0^\epsilon, u_0^\epsilon)$ defined on $\mathbb{T}$, such that $\rho_0^\epsilon(x)\geq \epsilon  \ \forall x, \ \|\rho_0-\rho_0^\epsilon\|_{L^1(\mathbb{T})}\rightarrow 0$, $\|u_{0}-u_0^\epsilon\|_\infty \rightarrow 0$ when $\epsilon\rightarrow 0$, then\\

 \textbf{Theorem 1.}\textit{ The solution of the system of ODEs (15-18), with $\alpha<\frac{1}{6}, 3\alpha+\beta<N-1$, with initial conditions $(\rho_0^\epsilon, u_0^\epsilon)$, provides a weak asymptotic solution (1) for the 1-D isothermal gas equations, which is global in time $t\in [0,+\infty[$ and in space $x\in \mathbb{T}$. }\\ 

 The result extends easily to 2-D and 3-D following section 6 in \cite{ColombeauODE}. We state the 2-D equations  (2, 3):

  $$\frac{d}{dt}\rho(x,y,t,\epsilon)=\frac{1}{\epsilon}[(\rho u^+)(x-\epsilon,y,t,\epsilon)-(\rho |u|)(x,y,t,\epsilon)+(\rho u^-)(x+\epsilon,y,t,\epsilon)+$$ \begin{equation}(\rho v^+)(x,y-\epsilon,t,\epsilon)-(\rho |v|)(x,y,t,\epsilon)+(\rho v^-)(x,y+\epsilon,t,\epsilon)]+ \epsilon^\beta,\end{equation}
for some $\beta>0$ to be made precise later,
$$\frac{d}{dt}(\rho u)(x,y,t,\epsilon)=\frac{1}{\epsilon}[(\rho u u^+)(x-\epsilon,y,t,\epsilon)-(\rho u |u|)(x,y,t,\epsilon)+(\rho u  u^-)(x+\epsilon,y,t,\epsilon)+ $$ \begin{equation}(\rho u v^+)(x,y-\epsilon,t,\epsilon)-(\rho u |v|)(x,y,t,\epsilon)+(\rho u  v^-)(x,y+\epsilon,t,\epsilon)] -\rho(x,y,t,\epsilon)\frac{\partial}{\partial x}\Phi(x,y,t,\epsilon),\end{equation}

 \textbf{ Numerical confirmations.}
 The discretization in time is standard from numerical schemes for ODEs: for fixed $\epsilon>0$ the Lipschitz properties of the ODEs (15-16) permit to prove convergence of the explicit Euler order one 
 method. The discretization in space is as follows: the space $\mathbb{R}^n$ is discretized
 by cells of length $\epsilon$ in each direction. The physical variables are constant in the cells. 
 Numerical tests from the explicit Euler order one method and the RK4 Runge Kutta method for the solution  of the ODEs (15, 16) have shown that the weak asymptotic method always gives the correct solutions, even  in presence of void regions, as this is the case in figures 1, 2 and 3.\\

 In figure 1 we present a demanding test from \cite{Bouchut} where an explicit solution is given. The numerical solution coincides with the explicit solution. 
 It is the Riemann problem $\rho_g=0, \rho_d=1, u_g=0= u_d$ at time $T=0.5$ and with $K=0.04$, which has been considered  in \cite{Bouchut} p. 154 and p. 157: this test is  difficult  due to the void region on the left in the initial conditions.
 The test has been done with  $dt=0.00002, \epsilon= 0.001,  2000$ space steps. \\

A striking fact is that one observes  that all the  technical ingredients $\beta, N, \Phi_\alpha, \rho_0^\epsilon>0$ in the construction of the weak asymptotic method of section 2 are really \  indispensible numerically: in the tests in \ figures \ 1 and 2 the

   \vskip 3.5 cm
     \begin{figure}[h]
 \includepdf[width=\textwidth]{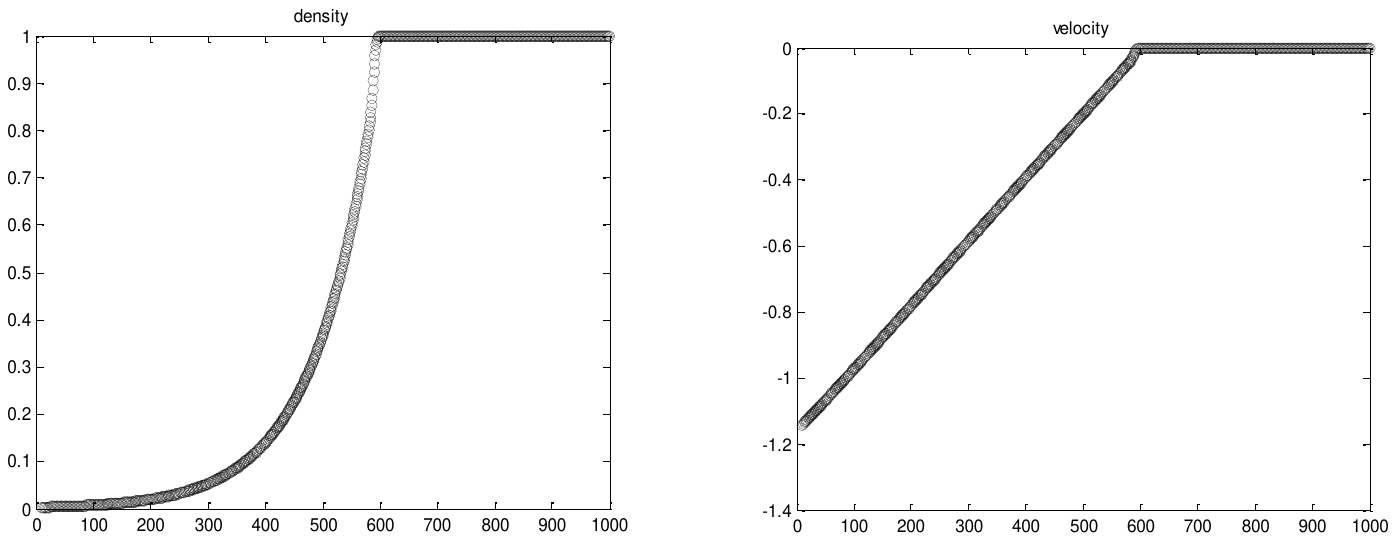}
\end{figure}
\textit{figure 1. A  demanding test on the system of isothermal gases with a void region; left: density, right: velocity.}\\
\\
numerical method fails to give a result if one of these ingredients is not taken into account. In figure 1 one has used $\rho_g=10^{-6}, \beta=10, N=2$ and an averaging on 3 cells for the state law (18) with coefficients 0.3, 0.4 and 0.3.  \\

 \textbf{3. Convergence to the analytic solution.}\\
 
In this section we present a slight modification of the weak asymptotic method in section 2, which has better properties, in particular one can prove that  when the initial data are analytic the weak asymptotic method  gives the classical analytic solutions.\\

 Let us consider a function $v(u), \ |u|<v(u)<|u|+\delta \ \forall u \in \mathbb{R},$ where $\delta>0$ is fixed, 
 and where the function $v$ is analytic, 
 for instance $v(u)=(u^2+\delta^2)^{\frac{1}{2}}$.
  Let us define $u^{\pm}$ by 
\begin{equation} u^+-u^-=u, \ \ u^++u^-=v(u) \end{equation}
instead of $ u^+-u^-=|u|$ in (13, 14). This gives
\begin{equation}u^+=\frac{v(u)+u}{2}=\frac{v(u)-|u|}{2}+\frac{|u|+u}{2}, \ u^-=\frac{v(u)-u}{2}=\frac{v(u)-|u|}{2}+\frac{|u|-u}{2}\end{equation}
i.e. $u^+$ and $u^-$ are larger by the quantity $\frac{v(u)-|u|}{2}$ than their  previous values considered in section 2. Now let us consider the ODE
\begin{equation} \frac{d}{dt}\rho(x,t,\epsilon)=\frac{1}{\epsilon}[(\rho u^+)(x-\epsilon,t,\epsilon)-(\rho v(u))(x,t,\epsilon)+(\rho u^-)(x+\epsilon,t,\epsilon)]+\epsilon^\beta\end{equation}
instead of (15), with a similar use of the new values $u^\pm$ and  $v(u)$  in (16) giving $ \frac{d}{dt}(\rho u)$.\\
\\
\\
\\
\\
\\
\\
\\
\\
\\
\\
\\
\\
\\
\\
\\
\\
\\
\\
\\
\\
\\
\\
\\
\\

 \begin{figure}[h]
\centering \includepdf[width=\textwidth]{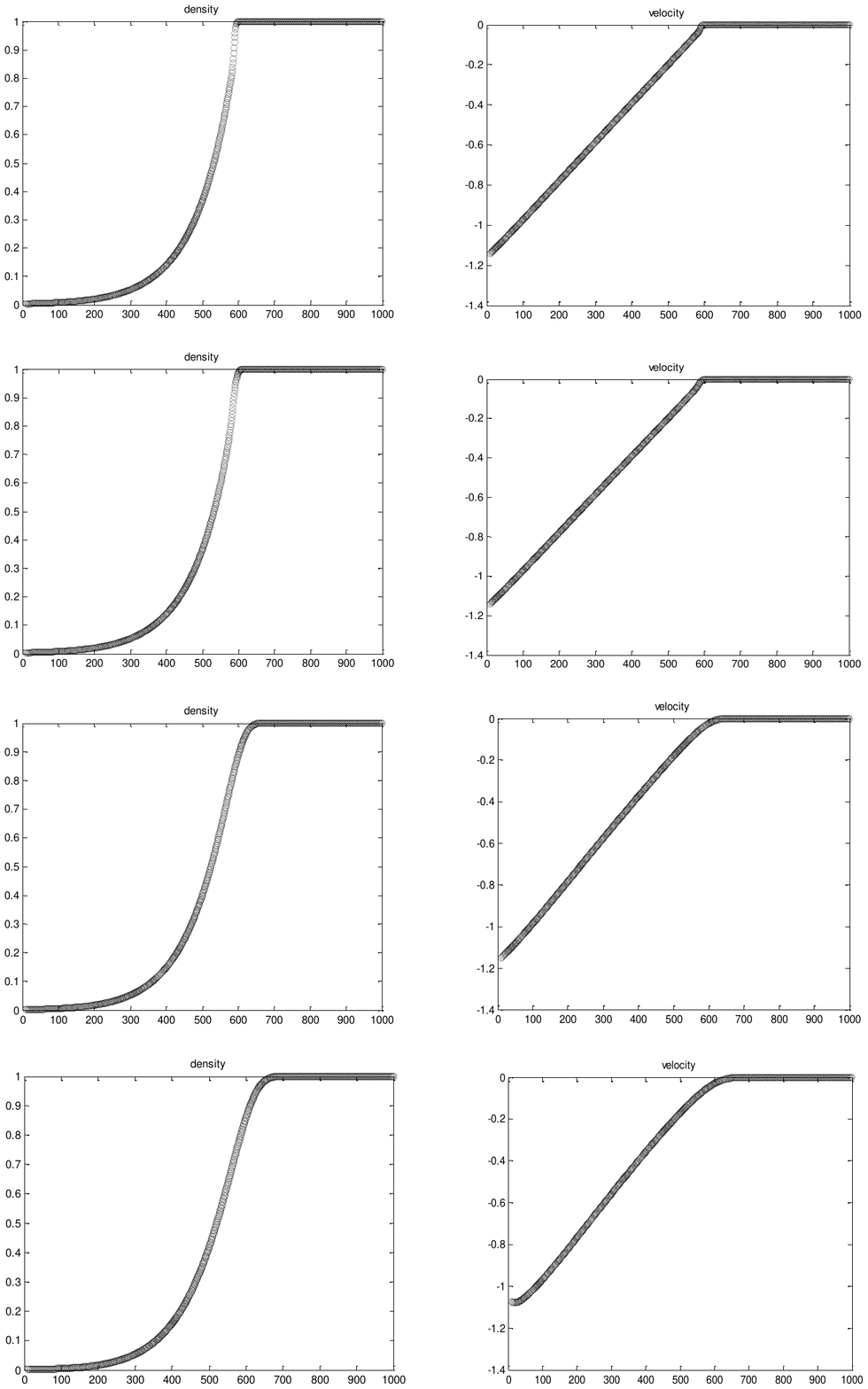}
\end{figure}
 \vskip 3 cm
\textit{figure 2. The test of figure 1 with variable values of the parameter $\delta=0.001,0.1,1,2$ from top to bottom.}\\
\\
 Let us consider a function $v(u), \ |u|<v(u)<|u|+\delta \ \forall u \in \mathbb{R},$ where $\delta>0$ is fixed, 
 and where the function $v$ is analytic, 
 for instance $v(u)=(u^2+\delta^2)^{\frac{1}{2}}$.
  Let us define $u^{\pm}$ by 
\begin{equation} u^+-u^-=u, \ \ u^++u^-=v(u) \end{equation}
instead of $ u^+-u^-=|u|$ in (13, 14). This gives
\begin{equation}u^+=\frac{v(u)+u}{2}=\frac{v(u)-|u|}{2}+\frac{|u|+u}{2}, \ u^-=\frac{v(u)-u}{2}=\frac{v(u)-|u|}{2}+\frac{|u|-u}{2}\end{equation}
i.e. $u^+$ and $u^-$ are larger by the quantity $\frac{v(u)-|u|}{2}$ than their  previous values considered in section 2. Now let us consider the ODE
\begin{equation} \frac{d}{dt}\rho(x,t,\epsilon)=\frac{1}{\epsilon}[(\rho u^+)(x-\epsilon,t,\epsilon)-(\rho v(u))(x,t,\epsilon)+(\rho u^-)(x+\epsilon,t,\epsilon)]+\epsilon^\beta\end{equation}
instead of (15), with a similar use of the new values $u^\pm$ and  $v(u)$  in (16) giving $ \frac{d}{dt}(\rho u)$.\\

First one can check that the proof of the  weak asymptotic method in section 2 holds without change: this is based  on the fact that the presence of any fixed value $\delta>0$ does not affect significantly the bounds and that one uses only the formula $u=u^+-u^-$ once one has replaced $v(u)$ by $u^++u^-$ in (48) and in the  ODE satisfied by $\rho u$. Second, numerical tests (figure 2) show that the presence of a fixed $\delta>0$   affects  the numerical results in the same way as a viscosity.\\

If we state $D(u)=\frac{v(u)-|u|}{2}$ and $\overline{u}^\pm=\frac{|u|\pm u}{2}$, i.e. $\overline{u}^\pm$ are the values $u^\pm$ used in section 2,  we obtain from (47, 48)\\

$ \frac{d}{dt}\rho(x,t,\epsilon)=\\
\\
\frac{1}{\epsilon}[(\rho \overline{u}^+)(x-\epsilon,t,\epsilon)-(\rho \overline{u}^+)(x,t,\epsilon)-(\rho \overline{u}^-)(x,t,\epsilon)+(\rho \overline{u}^-)(x+\epsilon,t,\epsilon)]+$\begin{equation}\epsilon^\beta+\frac{1}{\epsilon}[(\rho D(u)) (x-\epsilon,t,\epsilon)-2(\rho D(u)) (x,t,\epsilon)+(\rho D(u))(x+\epsilon,t,\epsilon)].\end{equation} 
The  first term in the second member gives the formula used in section 2. The last term is some kind of vanishing viscosity.\\
 
\textbf{Numerical confirmation.} In figure 2 the discretization is exactly the same as in section 2. The difference is the use of $v(u)=(u^2+\delta^2)^{\frac{1}{2}}$ instead of $|u|$. We
 observe that $\delta$  has the same influence as a viscosity coefficient: for $\delta=0.001$ (top panel) and  $\delta=0.1$ one obtains the exact solution; for $\delta=1$ one observes some viscosity effects which become quite important for $\delta=2$ (bottom panel).\\

Now from the replacement of $|u|$ by an analytic regularization $v(u)$ we prove that the weak asymptotic method gives the classical analytic solution at the limit $\epsilon\rightarrow 0$, which was  proved in \cite{ColombeauJDE} in a linear case when  the function $u$ has a fixed sign to avoid the singularity of the function absolute value at 0. The proof is given in the 1-D case since the multidimensional case is identical. The proof consists in applying the abstract nonlinear Cauchy-Kovalevska theorem of Nirenberg and Nishida \cite{Nishida} for each $\epsilon>0$ small enough, with  results uniform  in $\epsilon$.\\

\textbf{Recall of an abstract nonlinear Cauchy-Kovalevska theorem \cite{Nishida}.}  By definition  a scale of Banach spaces is a family of Banach spaces $(E_s)_s, 0<s\leq s_0$, such that $\forall s,s'\in ]0,s_0], \ s>s' \Rightarrow E_s\subset E_{s'}$ with inclusion of norm $\leq 1$. Let  $v_0\in E_{s_0}$ be given. If $R>0$ we denote by $B_s(v_0,R)$ the open ball in the Banach space $E_s$ of center $v_0$ and radius $R$.\\

Let $(E_s)_{0<s\leq s_0}$ be a scale of complex Banach spaces. Consider the Cauchy problem
 \begin{equation} \frac{dv}{dt}(t)=G(v(t)), \ \ t\in \mathbb{C}, \ v(0)=v_0.\end{equation}
Assume the existence of $R>0, C>0,K>0$ such that  properties i) and ii) hold\\
\\
i) $\forall s,s' \  /  \ 0<s'<s<s_0$ the map $u\longmapsto G(u)$ is holomorphic from $B_s(v_0,R)$ into $E_{s'}$, and satisfies a Lipschitz property in the sense

\begin{equation}\forall u\in B_s(v_0,R) \ \ \|G'(u)\|_{L(E_s,E_{s'})}\leq \frac{C}{s-s'}.
\end{equation}
\\
ii) $ \forall s<s_0 \ \  G(v_0)\in E_{s}$ and $\|G(v_0)\|_s\leq\frac{K}{s_0-s} $.\\

The abstract Cauchy-Kovalevska theorem states:
\textit{Then $ \exists \ a \ number \  a>0$ and a unique holomorphic function $t\longmapsto v(t)$ which $\forall s<s_0$  maps  $\{t\in\mathbb{C} / |t|<a(s_0-s)\}$ into $B_s(v_0,R)$ and is  solution of (50).}\\

The domain of $v$ depends on its range through the number $s$: one understands that the functions $v$ relative to various values of $s$ stick together. The proof is a holomorphic  form of the theorem in \cite{Nishida}, setting there $u=v-v_0, F(u,t)=G(v)$. Property (51) implies the Lipschitz property stated in \cite{Nishida} if $u,v \in B_s(v_0,R)$:
 \begin{equation}\|G(u)-G(v)\|_{s'}\leq \frac{C}{s-s'}\|u-v\|_s \ \forall u,v\in B_s(v_0,R).\end{equation} 
The method of proof is the  iteration method with  adequate  bounds, see \cite{Nishida}. The successive iterates lie in $B_s(v_0,R)$. The number $a>0$ depends only on $R, K,C $ (formulas 13 p. 630 and end of p. 632 in \cite{Nishida}).$\Box$\\

Now we can prove the following  coherence result.\\ 
\\
\textbf{Proposition 2.}\textit{ We consider  real valued analytic initial data $(\rho_0, \rho_0 u_0)$, independent of $\epsilon$, satisfying $\rho_0(x)>0 \ \forall x\in \mathbb{T}$. Then if $|t|$ is small enough the  solution of (15-18) with the modification  (46) tends to the classical analytic solution when  $\epsilon\rightarrow 0$.}\\ 

Proof. We denote by $\mathcal{H}(\mathbb{T}\times]-s,s[)$, respectively $\mathcal{C}(\mathbb{T}\times[-s,s])$, the spaces of all holomorphic, respectively continuous, functions $f=f(x,y), z=x+iy$ on the open strip $\mathbb{R}\times]-s,s[\subset \mathbb{C}$, respectively  the closed strip $\mathbb{R}\times[-s,s], x\in \mathbb{R}, |y|<s$ or $y\leq s$, periodic with period $2\pi$ in $x$-variable. For fixed $\epsilon$ we apply the abstract Cauchy-Kovalevska theorem above with the  classical scale of Banach spaces 

  \begin{equation} E_s=\{(f,g)\in (\mathcal{H}(\mathbb{T}\times]-s,s[)\cap \mathcal{C}(\mathbb{T}\times[-s,s]))^2\}\end{equation}  equipped with the  norm  \begin{equation}\|(f,g)\|_s=sup_{x\in \mathbb{T}, |y|\leq s } (|f(z)|, |g(z)|).\end{equation}  

The real number $s_0>0$ is chosen small enough so that the initial conditions $v_0^\epsilon:=v_0=(\rho_0, \rho_0 u_0)$, which are  independent on $\epsilon$ and are  holomorphic extensions of the given real analytic initial data, are elements of the space $E_{s_0}$. \\

Since we assume $\rho_0(x)>0 \ \forall x\in \mathbb{T}$ and since $\mathbb{T}$ is compact, one can choose $R>0$ small enough  so that if $(\rho,\rho u)\in B_s(v_0, R)$ then  $\forall s>0 \ small \ enough \  \exists \  b>0 \  /  \ |\rho(z)|>b \ \forall z=x+iy, x\in \mathbb{T}, |y|<s$, which permits  division by $\rho$ at each step in the iteration to obtain $u=\frac{\rho u}{\rho}$.\\

 We will also choose $R>0$ and $s_0>0$  small enough so that, 
 for a fixed $\delta>0$ in the definition of $v(u)=(u^2+\delta^2)^{\frac{1}{2}}$, $\forall (\rho,\rho u)\in B_s(v_0,R)$ the function $v$ is defined on the set $\frac{\rho u}{\rho}(\mathbb{T}\times]-s,s[)$ and bounded there, i.e. the set $\{(\frac{\rho u}{\rho}(x+iy))^2+\delta^2\}_{x\in \mathbb{T}, |y|<s<s_0}$ should remain at a $>0$ distance of the negative real demi-axis. This is possible since $\rho_0(x)$ and  $u_0(x)$ are real valued $\forall x\in \mathbb{T}$ and since  $\mathbb{T}$ is compact.\\

Then for any fixed $\epsilon>0$ small enough we consider the map $G_\epsilon$ defined from (15, 16) on $B_s(v_0, R)  \  \forall s<s_0$ with values in $E_{s'}, s'<s$ by \\

$G_\epsilon: [z\longmapsto(\rho(z), (\rho u)(z)]\in B_s(v_0, R)\longmapsto [[z\longmapsto\{ \frac{1}{\epsilon}[(\rho. u^+)(z-\epsilon)-\\
\\
(\rho.(u^++u^-))(z)+
(\rho. u^-)(z+\epsilon)], \ \frac{1}{\epsilon}[(\rho u. u^+)(z-\epsilon)-(\rho u.(u^++u^-))(z)+$\\

 \begin{equation}(\rho u. u^-)(z+\epsilon)]-\rho(z)\frac{\partial}{\partial x}\Phi(z)\}]]\in E_{s'}. \end{equation}

 For fixed $\epsilon$ we apply the abstract Cauchy-Kovalevska theorem with the analytic initial condition $v_0=(\rho_0,\rho_0 u_0)$ independent on $\epsilon$ and with the map $G_\epsilon$. We  obtain a solution $v^\epsilon=(\rho^\epsilon, (\rho u)^\epsilon)$   defined for $|t|<a(s_0-s)$ taking values in  $B_s(v_0,R) \ \forall s<s_0$ (the domain of $v^\epsilon$ depends of $s$: the same abuse of language has been done in the statement of the theorem above).
The key of the proof lies in that the assumptions of the abstract Cauchy-Kovalevska theorem are satisfied uniformly in $\epsilon$, therefore the properties (domains and bounds) of the solution will be also independent of $\epsilon$.  \\

The verification of the assumptions is based on  the mean value theorem and Cauchy's inequalities for holomorphic functions, \cite{Treves} p. 145. This can be checked easily: for simplification if $G_\epsilon(\rho, \rho u)$ would be  the function $z\longmapsto\frac{\rho u(z-\epsilon)-\rho u(z)}{\epsilon}$ then $\frac{\partial}{\partial(\rho u)} G_\epsilon(\rho,\rho u).w$ would be the function $z\longmapsto\frac{  w(z-\epsilon)- w(z)}{\epsilon}$ which is bounded as $( w)'(z-\theta(\epsilon,z) \epsilon), \ 0<\theta(\epsilon,z)<1$, therefore this amounts to a derivative
and  Cauchy's formula shows that $\frac{\partial}{\partial(\rho u)} G_\epsilon(\rho,\rho u)$ would map $E_s$ into $E_{s'}, \ s>s'$ with operator norm $\leq\frac{1}{s-s'}$ uniformly in $\epsilon$.\\

 The detailed formulas are  more complicated due to the presence of  $u=\frac{\rho u}{\rho}$ in $u^\pm$ inside the formula of $G_\epsilon$, and from the definition of $u^\pm$ in (46) through the function $v$. To clarify we give the formulas in the particular case $u\geq 0$ and $v(u)=|u|=u$ which implies  $u^+=u$ and $u^-=0$ and avoids longer formulas.  Then from (58)
 \begin{equation}G_\epsilon(\rho, \rho u)=z\longmapsto \{\frac{\rho u(z-\epsilon)-\rho u(z)}{\epsilon}+\epsilon^\beta,  \frac{\rho u^2(z-\epsilon)-\rho u^2(z)}{\epsilon}-\rho(z)\frac{\partial}{\partial x}\Phi(z)\}\end{equation} 
where $\Phi(z)=K[log(\rho(.)+\epsilon^N)*\phi_{\epsilon^\alpha}](z)$ from (18).\\

 To simplify the notation we set $X=\rho, Y=\rho u$; then
 \begin{equation}G_\epsilon(X, Y)=z\longmapsto \{\frac{Y(z-\epsilon)-Y(z)}{\epsilon}+\epsilon^\beta,  \frac{\frac{Y^2}{X}(z-\epsilon)-\frac{Y^2}{X}(z)}{\epsilon}-X(z)\frac{\partial}{\partial x}\Phi(z)\}.\end{equation} 
Therefore \\
$DG_\epsilon(X,Y).(w_1,w_2)=z\longmapsto
 \{\frac{w_2(z-\epsilon)-w_2(z)}{\epsilon}, \frac{  \frac{2Yw_2}{X}(z-\epsilon)-\frac{2Yw_2}{X}(z)    }{\epsilon }-\frac{   \frac{Y^2w_1}{X^2}(z-\epsilon)- \frac{Y^2w_1}{X^2}(z)   }{\epsilon}-$

 \begin{equation}Kw_1(z)[log(X+\epsilon^N)*(\phi_{\epsilon^\alpha})'](z)-KX(z)[\frac{w_1}{X+\epsilon^N}*(\phi_{\epsilon^\alpha})'](z)\}. \  \  \  \  \  \  \  \  \  \  \  \  \ \end{equation} 

If $(X,Y)\in B_s(v_0,R)\subset E_s$ then the values $X(x+iy), Y(x+iy)$ are defined for $x\in \mathbb{T}, |y|<s$ and are bounded in sup norm by $\|v_0\|_s+R$. Further $X(x+iy)>b>0$; the terms $\epsilon^N$ and $\epsilon^\beta$ in (60, 61) are not used to check the validity of (54) for $G_\epsilon$ uniformly in $\epsilon$.\\

The  abstract Cauchy-Kovalevska theorem asserts the existence of a solution $v^\epsilon: t\longmapsto (\rho^\epsilon(.,t), (\rho u)^\epsilon (.,t))\in B_s(v_0,R)$ to the ODEs (15, 16) defined on a domain of time  $\{|t|<a(s_0-s)\} $ which is independent of $\epsilon$ since domains and bounds on the data are uniform in $\epsilon$. Since the domains and bounds of the solutions are independent of $\epsilon$ the set of the functions $v^\epsilon$ is a normal family of holomorphic functions on $\{|t|<a(s_0-s)\}\times\mathbb{T}\times]-s,+s[$ for some  $a>0$ independent of $\epsilon$. From any sequence $(v^{\epsilon_n})_n$ one can extract a convergent subsequence. This subsequence converges to the classical analytic solution of the system of PDEs (10, 11). Therefore the whole family $(v^{\epsilon})$ converges to the classical analytic solution.$\Box$\\

    \textbf{4. Sequence of approximate solutions to the system of isentropic gas equations}.\\

      The  system of isentropic gas equations (2, 3, 4, 6) is stated in the form 
       \begin{equation} \frac{\partial}{\partial t}\rho+\frac{\partial}{\partial x}(\rho u)=0,  \end{equation} 
   \begin{equation}\frac{\partial}{\partial t} (\rho u)+\frac{\partial}{\partial x}(\rho u^2)+\rho \frac{\partial}{\partial x}\Phi=0,  \end{equation} 
   \begin{equation}\Phi=\frac{K\gamma}{\gamma-1} \rho^{\gamma-1}, K\geq 0 \ given, 1<\gamma\leq 2.   \end{equation}  
      
  We state the system of ODEs  as follows: we state the ODE (15) without the $\epsilon^\beta $ term which is no more needed because the state law (64) is different from (12), i.e.
  \begin{equation}\frac{d}{dt}\rho(x,t,\epsilon)=\frac{1}{\epsilon}[(\rho u^+)(x-\epsilon,t,\epsilon)-(\rho |u|)(x,t,\epsilon)+(\rho u^-)(x+\epsilon,t,\epsilon)],\end{equation}
(in which $u\pm$ are given by (13, 14), but could also be given by (46, 47) with the replacement of $|u|$ by $v(u)$).
  Then we state  the ODE (16), the formula  (17), and we replace the formula (18)  by   
    \begin{equation}  \Phi(x,t,\epsilon)= \frac{K\gamma}{\gamma-1}[(\rho(.,t,\epsilon))^{\gamma-1}*\phi_{\epsilon^\alpha}](x).\end{equation} 
     
    As usual the  convolution in (66) is justified by the fact that the state law is obtained from measurements which are always done in a space  region which is not too small. As in section 2 we assume  initial conditions $\rho_0^\epsilon$ and  $u_0^\epsilon$ defined on  $\mathbb{T}$. To obtain the a priori inequalities we assume (19-21). Then we obtain the same a priori inequalities as in proposition 1:     
    \textit{the statements (22 without the $\epsilon^\beta$ term), (23, 24, 25 with $2\alpha$ in place of $3\alpha$) and (26) hold as in proposition 1.}\\
    
 \textit{proof.} The proof of (22) in the present case is identical to the proof in proposition 1 without the $\epsilon^\beta$ term. For the proof of (23) here we have
 $$\frac{\partial}{\partial x}\Phi(x,t,\epsilon)=  \frac{K\gamma}{\gamma-1}\int[\rho(x-y,t,\epsilon)]^{\gamma-1} \frac{1}{\epsilon^{2\alpha}}\phi'(\frac{y}{\epsilon^\alpha})dy.$$
 Since $0<\gamma-1\leq 1$  and since  from (22) $\rho\in L^1(\mathbb{T})$ with $\|\rho\|_{L^1(\mathbb{T})}$ independent of $t$ and $\epsilon$, then a fortiori $\rho^{\gamma-1} \in L^1(\mathbb{T})$ with bounds independent of $t$ and $\epsilon$. Therefore
  $$  |\frac{\partial}{\partial x}\Phi(x,t,\epsilon)| \leq \frac{K\gamma}{\gamma-1}\|\rho^{\gamma-1}\|_{L^1(x-supp\phi)}\frac{1}{\epsilon^{2\alpha}}\|\phi'\|_\infty \leq \frac{const}{\epsilon^{2\alpha}}.$$
  The proofs of (24) and (26) in the present context are identical to those in proposition 1. $\Box$\\
  
  From the a priori estimates one obtains existence and uniqueness of a global solution in time $t\in [0,+\infty[$ and in space $x\in \mathbb{T}$. The weak asymptotic method for the analogs of (36, 37) is proved as in section 2; here it suffices to have $\alpha<\frac{1}{4}$ since (23) is stated  with $2\alpha$. For the state law (64) we have to prove that $\forall \psi \in \mathcal{C}_c^\infty(\mathbb{R}), \ \forall t$\\
  $$\int[\Phi(x,t,\epsilon)-\frac{K\gamma}{\gamma-1} (\rho(x,t,\epsilon))^{\gamma-1}]\psi(x)dx\rightarrow 0 $$  
  when $\epsilon\rightarrow 0$ where $\Phi(x,t,\epsilon)$ is given by (66). This  integral is equal to\\
  
  $\frac{K\gamma}{\gamma-1} \int_x[\int_y(\rho(x-y,t,\epsilon))^{\gamma-1}\frac{1}{\epsilon^\alpha}\phi(\frac{y}{ \epsilon^\alpha})dy -(\rho(x,t,\epsilon))^{\gamma-1} ]\psi(x)dx$=\\
\\
  $\frac{K\gamma}{\gamma-1} \int_{x,\mu}[(\rho(x-\epsilon^\alpha\mu,t,\epsilon))^{\gamma-1} -(\rho(x,t,\epsilon))^{\gamma-1} ]\phi(\mu)\psi(x)d\mu dx$=\\
\\
   $\frac{K\gamma}{\gamma-1} \int_{x,\mu}(\rho(x,t,\epsilon))^{\gamma-1} \phi(\mu)[\psi(x+\epsilon^\alpha \mu)-\psi(x)]d\mu dx=O(\epsilon^\alpha)$ \\
\\
from  the $L^1_{loc}(\mathbb{R})$ integrability of $\rho$ and since $0<\gamma-1\leq 1$ which permits integrability in $\rho^{\gamma-1}$. $\Box$\\ 
    
   Finally we have obtained: let the initial conditions $\rho_0\in L^1 (\mathbb{T})$, more generally a positive bounded Radon measure, $u^0\in L^\infty (\mathbb{T})$ and $\rho_0\geq 0$. Approximate the initial conditions by a family  $\rho_0^\epsilon, u_0^\epsilon \in \mathcal{C}(\mathbb{T}), \  \rho_0^\epsilon (x)>0  \  \forall x, \ \|\rho_0-\rho_0^\epsilon\|_{L^1(\mathbb{T})}\rightarrow 0$ and  $\|u_0-u_0^\epsilon\|_\infty \rightarrow 0$ when $\epsilon\rightarrow 0$. Then if $\alpha<\frac{1}{4}$:\\  
  
  \textbf{Theorem 2.}\textit{ The solution of the system of ODEs (65, 16, 17, 66)  with initial conditions $(\rho_0^\epsilon, u_0^\epsilon)$ provides a weak asymptotic solution (1) for the 1-D isentropic gas equations   (62-64), which is global in time $t\in [0,+\infty[$ and in space $x\in \mathbb{T}$.}\\

  \textbf{5. Sequence of approximate solutions to the system of selfgravitating collisionnal gases.} \\

Since the 1-D system is far simpler and the 3-D system quite analogous to the 2-D system, we state the system of isothermal collisionnal selfgravitating gases (\cite{Charru} p. 49, \cite{Coles} p. 207, \cite{Peacock} p. 460, \cite{Peter} p. 231) in 2-D for convenience:

 \begin{equation}\frac{\partial}{\partial t}\rho+\frac{\partial}{\partial x}(\rho u)+\frac{\partial}{\partial y}(\rho v)=0,   \end{equation}
  \begin{equation}\frac{\partial}{\partial t} (\rho u)+\frac{\partial}{\partial x}(\rho u^2)+\frac{\partial}{\partial y}(\rho uv)+\rho\frac{\partial}{\partial x}(\Phi_{press}+\Phi_{grav})=0,  \end{equation}
  \begin{equation}\frac{\partial}{\partial t} (\rho v)+\frac{\partial}{\partial x}(\rho uv)+\frac{\partial}{\partial y}(\rho v^2)+\rho\frac{\partial}{\partial y}(\Phi_{press}+\Phi_{grav})=0,  \end{equation}
  \begin{equation} \Phi_{press}=K log \rho, \ \ K\geq 0 \ constant,  \end{equation}
 \begin{equation} \Phi_{grav}(x,y,t)=const\int\rho(\xi, \eta,t) log\sqrt{(x-\xi)^2+(y-\eta)^2} d\xi d\eta.  \end{equation}
 where $\rho,u,v$ are as usual the density and the components of the velocity vector, $\Phi_{press}$ and $\Phi_{grav}$ are respectively the density of body force per unit mass caused  by the pressure and the gravitation. The system is similar in 3-D with (71) replaced by the usual 3-D elementary solution 
of the Poisson equation, $\frac{const}{r}, \ r=\sqrt{x^2+y^2}$, instead of $const. log r$ in 2-D.\\

The weak asymptotic solutions in 1-D, 2-D, 3-D for the system (67-71) are  adaptations of the  ones in section 2 in absence of selfgravitation, in which $\Phi(x,y,t,\epsilon)$ there is replaced by the sum $\Phi_{press}(x,y,t,\epsilon)+\Phi_{grav}(x,y,t,\epsilon)$ in which $\Phi_{press}$ is given by (18) and in which
\begin{equation} \Phi_{grav}(x,y,t,\epsilon)=[\Phi_{grav}(.,.,t)*\phi_{\epsilon^\alpha}](x,y).\end{equation}

The proofs (on the torus $\mathbb{T}^n, n=1,2,3$, always considered in this paper for simplification) take place in the Banach space $\mathcal{C}(\mathbb{T}^n), n=1,2,3$, and are not modified relatively to the isothermal case above, using the multidimensional pattern in section 6 of \cite{ColombeauODE}. One obtains again a weak asymptotic solution.\\

Remark. For the modelling of large scale structure formation of the universe one uses  periodic assumptions to model the cosmological principle that the universe is isotropic at large scales \cite{Coles} p. 305. Therefore, for this problem, the study on the torus $\mathbb{T}^3$ is more realistic than on the whole space $\mathbb{R}^3$.\\



 \textbf{6.  Sequence of approximate solutions to  the shallow water equations.}\\

 In this section 
    we construct a weak asymptotic solution for the 2-D shallow water equations stated in the form
 \begin{equation}\frac{\partial}{\partial t}h+\frac{\partial}{\partial x}(hu)+\frac{\partial}{\partial y}(hv)=0,   \end{equation}
  \begin{equation}\frac{\partial}{\partial t} (hu)+\frac{\partial}{\partial x}(hu^2)+\frac{\partial}{\partial y}(huv)+h\frac{\partial}{\partial x}\Phi=0,  \end{equation}
  \begin{equation}\frac{\partial}{\partial t} (hv)+\frac{\partial}{\partial x}(huv)+\frac{\partial}{\partial y}(hv^2)+h\frac{\partial}{\partial y}\Phi=0,  \end{equation}
 \begin{equation} \Phi=g(h+a),  \end{equation}
 where $h=h(x,y,t)$ is the water elevation, $(u,v)=(u(x,y,t),v(x,y,t))$ is the velocity vector in the $x,y$ directions respectively, $a=a(x,y)$ is the bottom elevation assumed to be of class $\mathcal{C}^2$ and $g=9.8$. The initial condition and the function $a$ are assumed to be periodic in the $x,y$ variables (with period $2\pi$ for convenience). We give a proof in 1-D since the 2-D extension is straightforward following the   pattern  exposed in section 6 of \cite{ColombeauODE}.\\
 
We state $h=\rho$ for convenience  because of the similarity with the other systems  considered before. Then the 1-D shallow water equations are:
   \begin{equation}\frac{\partial}{\partial t} \rho+\frac{\partial}{\partial x}(\rho u)=0,  \end{equation} 
   \begin{equation}\frac{\partial}{\partial t} (\rho u)+\frac{\partial}{\partial x}(\rho u^2)+\rho\frac{\partial}{\partial x} \Phi=0,  \end{equation} 
   \begin{equation}\Phi=g(\rho+a).   \end{equation} 
  We approximate (77-79) by the following system of ODEs:
\begin{equation}\frac{d}{dt}\rho(x,t,\epsilon)=\frac{1}{\epsilon}[(\rho u^+)(x-\epsilon,t,\epsilon)-(\rho |u|)(x,t,\epsilon)+(\rho u^-)(x+\epsilon,t,\epsilon)],\end{equation}

\begin{equation}\frac{d}{dt}(\rho u)(x,t,\epsilon)=\frac{1}{\epsilon}[(\rho u u^+)(x-\epsilon,t,\epsilon)-(\rho u |u|)(x,t,\epsilon)+(\rho u  u^-)(x+\epsilon,t,\epsilon)] -\rho(x,t,\epsilon)\frac{\partial}{\partial x}\Phi(x,t,\epsilon),\end{equation}
\begin{equation} \Phi(x,t,\epsilon)=g.(\rho(.,t,\epsilon)*\phi_{\epsilon^{\alpha}})(x)+g.a(x),\end{equation}
\begin{equation} u(x,t,\epsilon)=\frac{(\rho u)(x,t,\epsilon)}{\rho(x,t,\epsilon)},\end{equation}
for which we will prove that $\rho(x,t,\epsilon)>0,$ thus permitting division.\\

  We assume $\rho_0$ and $u_0$ are given  with the properties $\rho_0 \in L^1(\mathbb{T})$ and  $u_0\in   L^\infty(\mathbb{T})$ and that $\rho_0^\epsilon$ and $u_0^\epsilon$ are  continuous regularizations of $\rho_0$ and $u_0$ respectively, with  respective $L^1$ and $L^\infty$ bounds  independent on $\epsilon$, and $\rho_0^\epsilon (x)>0 \  \forall x\in \mathbb{T}$.\\

 We first establish a priori inequalities to prove existence of a global solution to (80-83).  For fixed $\epsilon>0$ we assume (19, 20, 21). Then one obtains the a priori inequalities:\\

  \textit{the statements (22 without the $\epsilon^\beta$ term), (23, 24, 25 with $2\alpha$ in place of $3\alpha$) and (26) hold as in proposition 1.}\\
  

      The proof of the analog of (22) is identical to the proof of (22) in proposition 1. For the proof of the analog of (23)     $$\frac{\partial}{\partial x}\Phi(x,t,\epsilon)=g[\rho(.,t,\epsilon)*(\phi_{\epsilon^\alpha})'](x)+ga'(x)$$
     i.e.
     $$\frac{\partial}{\partial x}\Phi(x,t,\epsilon)=g\int\rho(x-y,t,\epsilon)\frac{1}{\epsilon^{2\alpha}}\phi'(\frac{y}{\epsilon^\alpha})dy+ga'(x).$$
     Therefore
     $$|\frac{\partial}{\partial x}\Phi(x,t,\epsilon)|\leq \frac{const}{\epsilon^{2\alpha}}\|\rho(.,t,\epsilon)\|_{L^1(x-supp(\phi))} +const.$$ $\Box$.\\ 
    The proofs of the analogs of (24, 26) are similar to those  in proposition 1.\\

 The existence of a global solution  for fixed $\epsilon$ is obtained from these a priori estimates.\\
 
 It remains to prove that the solution of the system of ODEs (80-83) provides a weak asymptotic solution for system (77-79) when $\epsilon\rightarrow 0$. For (77, 78) the proof in the isothermal case applies. For the state law (79) we have to check that
 
\begin{equation}  \int\Phi(x,t,\epsilon)\psi(x)dx=g\int\rho(x,t,\epsilon)\psi(x)dx+g\int a(x)\psi(x)dx+f(\epsilon), \end{equation}
where $\Phi$ is given by (82) and 
where $f(\epsilon)\rightarrow 0$ when $\epsilon\rightarrow 0$.
 To check (84) it suffices from (82) to consider the integral \\

$\int[\rho(.,t,\epsilon)*\phi_{\epsilon^\alpha}](x)\psi(x)dx= \int \rho(x-y,t,\epsilon)\frac{1}{\epsilon^\alpha}\phi(\frac{y}{\epsilon^\alpha})\psi(x)dydx=\int \rho(x,t,\epsilon)\phi(\mu)$\\
\\
$\psi(x+\epsilon^\alpha\mu)d\mu dx=\int \rho(x,t,\epsilon)\psi(x)dx+\int \rho(x,t,\epsilon)\phi(\mu)[\psi(x+\epsilon^\alpha \mu)-\psi(x)]d\mu dx.$ \\
\\
The last term is $O(\epsilon^\alpha)$ from the $L^1_{loc}$ property (22) of $\rho$  since integration takes place on a compact set. $\Box$ \\

  Finally we have proved: let the initial conditions $\rho_0\in L^1 (\mathbb{T}), u^0\in L^\infty (\mathbb{T}), \  \rho_0\geq 0,  \ a$ a function of class  $\mathcal{C}^2$ on $\mathbb{T}$. Approximate the initial conditions by a family $(\rho_0^\epsilon, u_0^\epsilon) \in (\mathcal{C}(\mathbb{T}))^2$, such that $\rho_0^\epsilon(x)>0  \ \forall x, \ \|\rho_0-\rho_0^\epsilon\|_{L^1(\mathbb{T})}\rightarrow 0$, $\|u_0-u_0^\epsilon\|_\infty \rightarrow 0$ when $\epsilon\rightarrow 0$. Then\\
  
  \textbf{Theorem 3}. \textit{The solution of the system of ODEs (80-83),  $0<\alpha<\frac{1}{4}$,  with initial conditions $(\rho_0^\epsilon, u_0^\epsilon)$ provides a weak asymptotic method (1) for the 1-D shallow water equations (77-79), which is global in time $t\in [0,+\infty[$ and in space $x\in \mathbb{T}$.}\\ 

The method extends to 2-D following \cite{ColombeauODE} section 6.   One obtains coherence with the classical analytic solution as in proposition 2. Numerical tests from \cite{Toro} have given results similar to those in section 2 and 3, with the possible use of $v(u)$ in place of $u$. In figure 3 we have tested the three Riemann problem tests 1, 2 and 3 (from top to bottom) in \cite{Toro} pp. 109-124. Top: we obtain the exact solution in test 1 with 500 space steps, dt=0.0001, an averaging on 3 cells to represent the convolution (82) on the state law with coefficients $\alpha, 1-2\alpha, \alpha$ with $\alpha=0.1$, no averaging needed on the initial condition and $\delta=0$ in the formula $v(u)=(u^2+\delta^2)^\frac{1}{2}$.
Middle: we obtain the exact solution in test 2 with 2000 space steps, dt=0.00004, an averaging for the state law (82) with $\alpha=0.1$, a similar  averaging with   $\alpha=0.1$ for the initial condition and $\delta=1$. Bottom: we obtain the exact solution in test 3 with 5000 space steps, $dt=8.10^{-6}$, averagings with $\alpha=0.1$ for the state law (82) and the initial condition, $\delta=0.5$. Due to a void region in initial condition ($\rho=1$ on the left and 0 on the right) one has replaced the right value $\rho=0$ by $10^{-10}$ since  the proof of weak asymptotic method requests $\rho_{0,\epsilon}(x)>0 \ \forall x$. \\

 \textbf{7. Conclusion.} \\

In view of providing a substitute of solutions that could explain the observed numerical results such as \cite{Lax1, Lax2}, we have constructed  approximate solutions (up to any given accuracy in the sense of distributions)  to the general Cauchy problem on the n-dimensional torus $\mathbb{T}^n, n=1,2,3$,    for  some standard equations of fluid dynamics in presence of shocks and void regions.\\

   Up to our knowledge, sequences of global approximate solutions to the general initial value problem for the equations of fluid dynamics with a complete mathematical proof  that they tend to satisfy the equations had not been constructed previously. The statement of the  method is based on a family of two nonlinear ordinary differential equations in a classical Banach space of functions, for which a priori estimates permit to prove existence-uniqueness of  global solutions (in space and in positive time).
The continuity equation (mass conservation) gives a $L^1$ control on density which, from its use inside the state law giving pressure, permits from the Euler equations some control on velocity.\\

The method of proof allows  numerical calculations from standard convergent numerical methods for ODEs such as the explicit Euler order 1 method and the RK4 method. They  have given  the known solutions in all   tests.\\
\\
\\
\\
\\
\\
\\
\\
\\
\\
\\
\\
\\
\\
 \vskip 7 cm
 \begin{figure}[h]
\centering \includepdf[width=\textwidth]{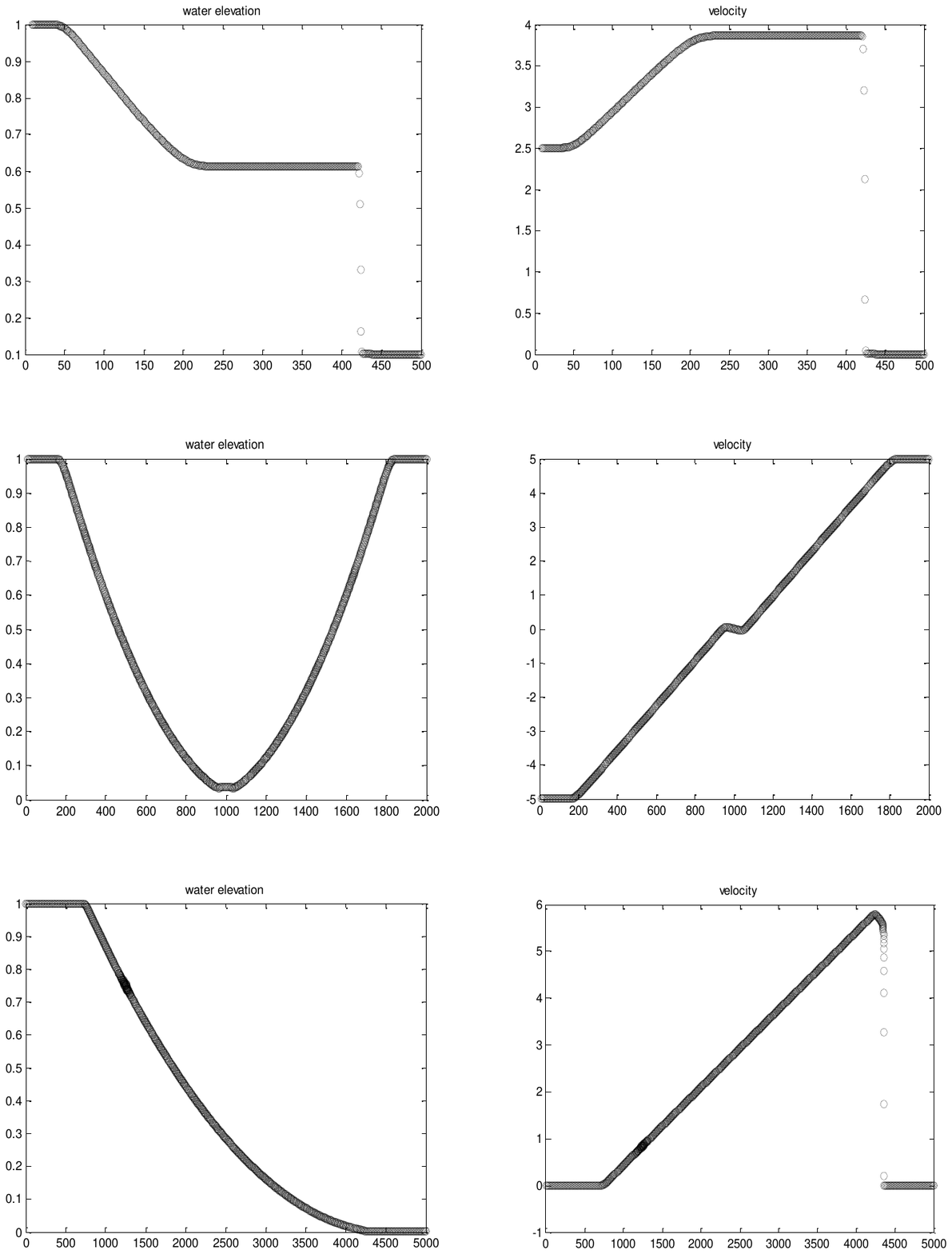}
\end{figure}

\textit{figure 3. From top to bottom tests 1, 2 and 3 on the equations shallow water in \cite{Toro} pp. 109-124; left: water elevation, right water velocity. The results were obtained with the explicit Euler order one method applied to the ODEs (80, 82). One observes coincidence with the exact solutions given in \cite{Toro}, although these tests have been chosen on purpose to be very demanding from the numerical viewpoint due to the presence of void regions in the water elevation.}\\

  Uniqueness of an "admissible" class of such sequences of approximate solutions, containing those constructed in this paper ("existence part"), such that all sequences in this class  give same values for all time at the limit $\epsilon\rightarrow 0$ ("uniqueness part") has been obtained so far only for some linear equations \cite{ColombeauCK, CCNO}. Passage to the limit by some kind of weak compactness has not been considered so far.\\
\\ Acknowledgements. The author is very grateful to members of the Instituto de Matematica of The Universidade de S\~ao Paulo, of the Instituto de Matematica, Estatistica e Computa\c c\~ao Cientifica of the Universidade Estadual de Campinas and of the Instituto de Matematica of the Universidade Federal do Rio de Janeiro for their attention, encouragements and suggestions while doing this work.\\

\textbf{Appendix}: \textit{Why one could be forced to use sequences of approximate solutions as a substitute of exact solutions.}\\

We recall known examples which suggest the possible absence of distribution  solutions also in some instances of fluid dynamics.\\

$\bullet$ The system 
\begin{equation} \frac{\partial u}{\partial t}+ \frac{\partial u^2}{\partial x}=0,\end{equation}
\begin{equation} \frac{\partial v}{\partial t}+ \frac{\partial (2 uv)}{\partial x}=0,\end{equation}
\begin{equation} \frac{\partial w}{\partial t}+ \frac{\partial (2 uw)}{\partial x}+ \frac{\partial v^{n+1}}{\partial x}=0,\end{equation}
produces shock waves involving a power $\delta^n$ of the Dirac delta measure as a continuation of analytic solutions after their blow up. To prove this fact we produce a sequence of approximate solutions with the same properties as those in this paper and we compute them from a convergent numerical scheme for ODEs. One observes numerically the powers $\delta^n$ \cite{ColombeauJDE}. System (85-86) is not strictly hyperbolic but this fact has only served to create a $\delta$-wave in $v$ to insert into (87) and such $\delta$-wave can be created as well in strictly hyperbolic systems: see $v$ in system (88-89) below. \\ 

$\bullet$ The Keyfitz-Kranzer system 
\begin{equation} \frac{\partial u}{\partial t}+ \frac{\partial (u^2-v)}{\partial x}=0,\end{equation}
\begin{equation} \frac{\partial v}{\partial t}+ \frac{\partial (\frac{u^3}{3}-u)}{\partial x}=0,\end{equation}
shows a more subtle phenomenon. The approximate solutions $u^\epsilon, v^\epsilon$ explicitely calculated in  formula (8) in \cite{Keyfitz} from the Dafermos-Di Perna viscosity, and obtained again numerically  in \cite{Sanders} from the usual viscosity  have a limit $u,v$ in the sense of distributions when $\epsilon\rightarrow 0$: $u$ shows a simple discontinuity and $v$ shows a Dirac mass located on the discontinuity; they move with constant speed. When one inserts such $u$ and $v$ in $u^2-v$ the Dirac measure in $v$ subsists untouched in $u^2-v$ and therefore the term $\frac{\partial (u^2-v)}{\partial x}$ has a nonzero $\delta'$ term. This term cannot be compensated in $ \frac{\partial u}{\partial t}$ because $u$ is a simple discontinuity. Equation (88) cannot be satisfied. Similarly in (89) $\frac{\partial v}{\partial t}$ shows a $\delta'$ term that cannot be compensated in $\frac{\partial (\frac{u^3}{3}-u)}{\partial x}$. At the level of the approximate sequence  $(u^\epsilon, v^\epsilon)$ with $\epsilon\not=0$ the compensations take place because the function $u^\epsilon$ has a term (the last term in the first formula (8) in \cite{Keyfitz}) which could be intuitively refered to has a "square root" of a Delta measure, i.e. an object whose square tends in the sense of distributions to the Dirac measure when $\epsilon \rightarrow 0$: then $u^2$ has a $\delta$-term that compensates the $\delta$ term of $v$ in $u^2-v$; a similar verification holds for equation (89). This shows that when one  passes to the limit $\epsilon\rightarrow 0$ in the approximate solutions in the sense of distributions  one obtains a result $(u,v)$ that, when plugged into the equations, does not always satisfy them.\\

 \end{document}